\documentclass[11pt,twocolumn ]{IEEEtran}
\usepackage{graphicx, color,  amssymb, amsmath, cite}
\usepackage{amsthm}
\usepackage[english]{babel}
\usepackage{cuted}
\usepackage{framed}

\makeatletter
\newcommand{\subjclass}[2][1991]{%
  \let\@oldtitle\@title%
  \gdef\@title{\@oldtitle\footnotetext{#1 \emph{Mathematics subject classification.} #2}}%
}

\makeatother

\newtheorem{thrm}{Theorem}[section]

\newtheorem{lem}[thrm]{Lemma}

\title{Computationally Efficient Bounds for the Sum of Catalan Numbers}

\author{ Kevin Topley \\  Email: kevint@ece.ubc.ca}

\date{}
\begin{document}

\maketitle

\begin{abstract} Easily computable lower and upper bounds are found for the sum of Catalan numbers. The lower bound is proven to be tighter than the upper bound, which previously was declared to be only an asymptotic. The average of these bounds is proven to be also an upper bound, and empirically it is shown that the average is superior to the previous upper bound by a factor greater than $(9/2)$. \end{abstract}


\section{Introduction}\label{sec:intro}

The Catalan numbers form a sequence of natural numbers that occur in a variety of counting problems  \cite{stan,cor90}. The sum of the first $n$ Catalan numbers has been shown to equal the number of paths starting from the root in all ordered trees with $(n+1)$ edges \cite{oeis2}. The sum of the first $n$ Catalan numbers also equals (a) the sum of the mean maximal pyramid size over all Dyck $(n+1)$-paths, and (b) the sum of the mean maximal saw-tooth size over all Dyck $(n+1)$-paths \cite{oeis3}.

Although there are numerous closed-form expressions for the $k^{th}$ Catalan number $C_k$, none of them are especially attractive from a computational standpoint \cite{dut86}. Determining the sum of the first $n$ Catalan numbers requires computation of $C_1, C_2, \ldots, C_n$ \cite{oeis4}; it is thus reasonable to search for an accurate and easily computable approximation to the sum of Catalan numbers. Motivated by its applications and cumbersome expression, we will find computationally efficient upper and lower bounds to the sum of Catalan numbers. The tightness of these approximations is quantified both analytically and empirically.

\section{Main Results}\label{sec:main}

The $k^{th}$ Catalan number $C_k$ is defined as, \begin{equation}\label{cat1} C_k = \frac{{2k \choose k}}{k+1}   = \prod_{i = 0}^{k-2} \Big( \frac{2k-i}{k-i} \Big) \ ,  \ k \in \{ 1, 2, \ldots \} \ . \end{equation} The sum of the first $n$ Catalan numbers is then given by $S_n$, \begin{equation}\label{cat2} S_n =  \sum_{k = 1}^{n}  C_k \ . \end{equation} The following asymptotic limit of $S_n$ has been proposed \cite{oeis},  \begin{equation}\label{cat3} S_n   \sim  \frac{4^{n+1}}{3 \sqrt{\pi n^3 }} \doteq u(n) \ . \end{equation} We will prove $(\ref{cat3})$ is actually an upper bound for $S_n$. Furthermore, we find a more accurate approximation to $S_n$ is given by the following lower bound, \begin{equation}\label{cat4} S_n >   \frac{4^{n+1}}{3(n+1) \sqrt{\pi  n }  }  \doteq  \vartheta(n) \ . \end{equation} Specifically, we have the result, \begin{equation}\label{main} \begin{array}{llll} &   \  u(n)  > S_n   > \vartheta(n) \ \ , \ \forall \ n \geq 1  \\ &   u(n) + \vartheta(n)   >  2  S_n \  \  , \  \forall \  n \geq 8 \ . \end{array} \end{equation} Note that the quotient $u(n)/\vartheta(n)$ approaches $1$ as $n$ approaches infinity, thus both approximations are asymptotically equal to $S_n$.

\subsection{Proof of Main Results}\label{sec:proof}

The main result $(\ref{main})$ is proven here in Thm.$\ref{thm1},\ref{thm2}$. To obtain $(\ref{main})$ we require Lemmas $\ref{lem1} - \ref{lem3}$. Let $\mathbb{N}$ to be the set of non-negative integers, and $\mathbb{R}$ the set of real numbers.

\begin{lem}\label{lem1} The sum of the first n Catalan numbers, $S_n$ (cf. $(\ref{cat2})$) has a lower bound $\ell_n \doteq 4 C_n /3 $. \end{lem}

\noindent  \emph{Proof of Lem.$\ref{lem1}.$} Rearranging $S_n > \ell_n$ yields the inequality, \begin{equation}\label{eq1a} S_n < 4 S_{n-1} \ . \end{equation} We will use the recurrence,  \begin{equation}\label{com1} C_k = \frac{2(2k -1)}{k+1} C_{k-1} \end{equation} which can be easily obtained from $(\ref{cat1})$ \cite{dut86}. Applying $(\ref{com1})$ to the Catalan numbers $C_k$ yields, $$ 4 C_{k-1} - C_k  = \frac{3 C_k  }{2k -1} > 0 \ , \ \forall  \ k \in \{ 1, 2, \ldots \} $$   thus we obtain $(\ref{eq1a})$. $\ \ \ \quad \quad \quad \quad \quad \quad \quad \quad \quad \quad \quad \quad \quad \quad  \ \ \quad \blacksquare$

\begin{lem}\label{lem2}  The sum of the first n Catalan numbers, $S_n$ (cf. $(\ref{cat2})$) has an upper bound $u_n$ (cf.$(\ref{cat3})$).   \end{lem}

\noindent   \emph{Proof of Lem.$\ref{lem2}.$} It is shown in \cite{dut86} that the $k^{th}$ Catalan number $C_k$ has an upper bound, \begin{equation} \label{ub} C_k < \frac{ 2^{2k+1}}{(k+1) \sqrt{ \pi (4k+1) }}  \doteq \nu(k) \ .  \end{equation} Numerical evaluation shows that for $n \in \{ 1, 2, \ldots, 12 \}$, $$ S_n < \sum_{k = 1}^n  \frac{ 2^{2k+1}}{(k+1) \sqrt{ \pi (4k+1) }} < u_n $$ where $(\ref{ub})$ has been applied to $S_n$. We now proceed by proving, \begin{equation}\label{mainlem} u_n + \ell_n > 2 S_n  \  , \ \forall  \ n \geq 13 \ . \end{equation} Subtracting $\ell_n$ from $(\ref{mainlem})$ and multiplying by $3$ yields, \begin{equation}\label{y1}   2 S_n + 4 S_{n-1}  < 3 u_n \ . \end{equation}  Numerical evaluation verifies $(\ref{y1})$ for $n = 13$.  Next we take $(\ref{y1})$ as the inductive assumption. It remains to be shown that,  \begin{equation}\label{y2}  3 u_{n+1} - 2 S_{n+1} -  4 S_n  > 0 \ .  \end{equation} Applying the inductive assumption $(\ref{y1})$ to $(\ref{y2})$ provides the sufficient condition, \begin{equation}\label{p1} 3 u_{n+1}  >  3 u_n + 2 C_{n+1} + 4 C_n \ . \end{equation}

\noindent  Applying $(\ref{com1})$ and $(\ref{ub})$ to the sum $2 C_{n+1} + 4 C_n$ we obtain the upper bound, \begin{equation}\label{array} \begin{array}{llll}  & 2C_{n+1} + 4 C_n = 4 C_n \big( 1 + \frac{2n+1}{n+2} \big) \\ & \ \ \ \ \ \ \ \ \ \ \ \ \ \ \ \  \ \ <    \frac{ 3 \cdot 2^{2n+3} }{ (n+2) \sqrt{ \pi (4n + 1) }} \ . \end{array} \end{equation} Applying $(\ref{array})$ to $(\ref{p1})$ and simplifying yeilds the sufficient condition, \begin{equation}\label{p3} 4 \geq \sqrt{  \Big( \frac{n+1}{n} \Big) ^3 } + \sqrt{ \frac{36 (n+1)^3}{(n+2)^2 (4n + 1) }  } \ . \end{equation} By expanding $(\ref{p3})$ we have, \begin{equation}\label{eq1}  ( \textbf{h} ' \textbf{n}  ) ^2 \geq 4 ( \textbf{q} ' \textbf{n}  \textbf{r} ' \textbf{n} )   \end{equation} where the $\mathbb{R}^{7 \times 1}$ vectors $\textbf{h}$,$\textbf{q}$,$\textbf{r}$, and $\textbf{n}$, are defined, $$ \begin{array} {llll}  & \textbf{h}  =  [ -4, -24, -84, -91 , 129 , 135, 24] ' \\ & \textbf{q} =  [ 0, 0, 0, 36 , 108 , 108, 36] ' \\  & \textbf{r} =  [ 4, 24, 84, 119 , 83 , 29, 4] ' \\ & \textbf{n} = [1, n, n^2, n^3, n^4, n^5, n^6] '  \ . \end{array} $$ Expanding and then simplifying $(\ref{eq1})$ yields, $$  \begin{array} {llll}  &  \quad  \quad  \quad  \quad  \quad  \quad \textbf{j}' \textbf{N} \geq 0   \\ & \ \   \textbf{j}  = [ 16 ,192 , 1248 , 4184,    5208,    \\ &  \ \ \ \ \ \ \  -16176 , -84431, -150414  ,\\ & \quad \quad \ \  \ -115497,   -35634 , -1791, 576 ] '  \\  & \ \ \ \ \ \  \textbf{N}_i = n^{i-1} \ , \ i \in \{1,2, \ldots, 12 \}  \end{array} $$ where $\textbf{N} \in \mathbb{R}^{12 \times 1}$ has $i^{th}$ element $\textbf{N}_i$. Denoting the $i^{th}$ element of $\textbf{j} \in \mathbb{R}^{12 \times 1}$ as $\textbf{j}_i$, it is clear that $\textbf{j}' \textbf{N} >0$ for $n = 13$ since $\textbf{j}_i >0$ for $i \in \{1,2, \ldots, 5\}$, $\textbf{j}_i < 0$ for $i \in \{6,7, \ldots, 11\}$,  and by numerical evaluation we have, $$ 0 < -  \sum_{i = 6}^{11} \big(  \textbf{j}_i / 13^{12-i} \big) < \textbf{j}_{12} \ . $$

We have shown $(\ref{mainlem})$ holds for all integers greater than $12$. Lemma $\ref{lem1}$ proves $\ell_n < S_n$, thus if $u_n \leq S_n$  then $\ell_n + u_n < 2 S_n$, which contradicts $(\ref{mainlem})$. $\ \quad \quad \quad \quad \ \quad \quad \quad \ \ \ \blacksquare$

\begin{lem}\label{lem3} The sum of the first n Catalan numbers, $S_n$ (cf. $(\ref{cat2})$) has a lower bound $\vartheta_n$ (cf.$(\ref{cat4})$). \end{lem}

\noindent  \emph{Proof of Lem.$\ref{lem3}.$} Subtracting $\ell_n$ from $\vartheta_n < S_n$ and multiplying by $3$ yields, \begin{equation}\label{a1}  3 \vartheta_n -  3 \ell_n  < 4 S_{n-1}  - S_n \ .  \end{equation} Rearranging $(\ref{a1})$ and applying Lem.$\ref{lem1},\ref{lem2}$ yields the sufficient condition, \begin{equation}\label{n1}  3 \vartheta_n \leq 4 \ell_{n-1} + 3 \ell_n  - u_n \ . \end{equation} \noindent It is shown in \cite{dut86} that the $k^{th}$ Catalan number $C_k$ has a lower bound, \begin{equation} \label{lb} C_k > \frac{ 2^{2k-1}}{k (k+1) \sqrt{ \pi / (4k-1) }} \ .  \end{equation} Applying $(\ref{com1})$ and $(\ref{lb})$ to the sum $4 \ell_{n-1} + 3 \ell_n$  we obtain, \begin{equation}\label{suf2}\begin{array}{llll}  & 4 \ell_{n-1} + 3 \ell_n  = 4 C_n  \Big( \frac{8n -1}{6n-3} \Big) \\ & \  \ \ \ \   >  \frac{2^{2n+1}}{n(n+1) \sqrt{ \pi / (4n-1) } } \Big( \frac{8n -1}{6n-3} \Big) \ . \end{array} \end{equation} Substituting $(\ref{suf2})$ as well as the expressions for $u_n$ (cf.$(\ref{cat3})$) and $\vartheta_n$ (cf.$(\ref{cat4})$) in $(\ref{n1})$ yields the sufficient condition, \begin{equation}\label{poly0} \begin{array}{llll} &  \frac{4^{n+1}}{ (n+1) \sqrt{ \pi n } }  \leq \\ &  \ \ \ \   \frac{2^{2n+1}}{n(n+1) \sqrt{ \pi / (4n-1) } } \Big( \frac{8n -1}{6n-3} \Big)  - \frac{4^{n+1}}{ 3 \sqrt{ \pi n^3 }} \ . \end{array}  \end{equation} Simplifying, $(\ref{poly0})$ becomes,  \begin{equation}\label{poly1} \begin{array}{llll} & 4 ( 4n+1)^2 (6 n - 3)^2 \\ & \ \ \ \ \ \ \ <  9 n (8n-1)^2 (4n-1) \ . \end{array}  \end{equation} Rearranging $(\ref{poly1})$ yields the equivalent condition $68 n^2 >  17 n + 4$, which holds for all $n \geq 1.$ $\quad \quad  \quad \quad  \quad \quad  \quad \quad \ \ \ \ \ \ \blacksquare$

\begin{thrm}\label{thm1} For all $n \geq 1$, $u_n  > S_n > \vartheta_n $.\end{thrm}

\noindent \emph{Proof of Thm.$\ref{thm1}$.} The result is a combination of Lem.$\ref{lem2}, \ref{lem3}$.  $ \   \quad \quad  \quad \quad  \quad \quad \quad \quad \quad \quad \quad \quad \quad \quad \quad \ \ \ \ \ \ \ \  \blacksquare$

\begin{thrm}\label{thm2} For all $n \geq 8$, $ u_n + \vartheta_n   >  2  S_n  $.\end{thrm} 

\noindent \emph{Proof of Thm.$\ref{thm2}$.} Numerical evaluation shows that for $n \in \{8,9, \ldots, 12 \}$,  $$ 2 S_n < 2 \sum_{k = 1}^n  \frac{ 2^{2k+1}}{(k+1) \sqrt{ \pi (4k+1) }} <  u_n + \vartheta_n  $$ where $(\ref{ub})$ has been applied to $S_n$. We now show that $\vartheta_n > \ell_n$, thus the result follows from $(\ref{mainlem})$. Apply $(\ref{ub})$ to $\ell_n$, $$ \ell_n <  \frac{ 2^{2n+3}}{3(n+1) \sqrt{ \pi (4n+1) }} \ . $$ It suffices to then show,  $$ \frac{ 2^{2n+3}}{3(n+1) \sqrt{ \pi (4n+1) }} <  \frac{4^{n+1}}{3(n+1) \sqrt{\pi  n }  }  $$ which simplifies to $1 > 0.$ $\quad \quad \quad \quad \quad \quad \quad \quad \quad  \quad \quad \quad \quad \ \  \blacksquare$

\section{Numerical Results}\label{sec:num}

\begin{figure}[htb] \center
\includegraphics[width=1 \linewidth]{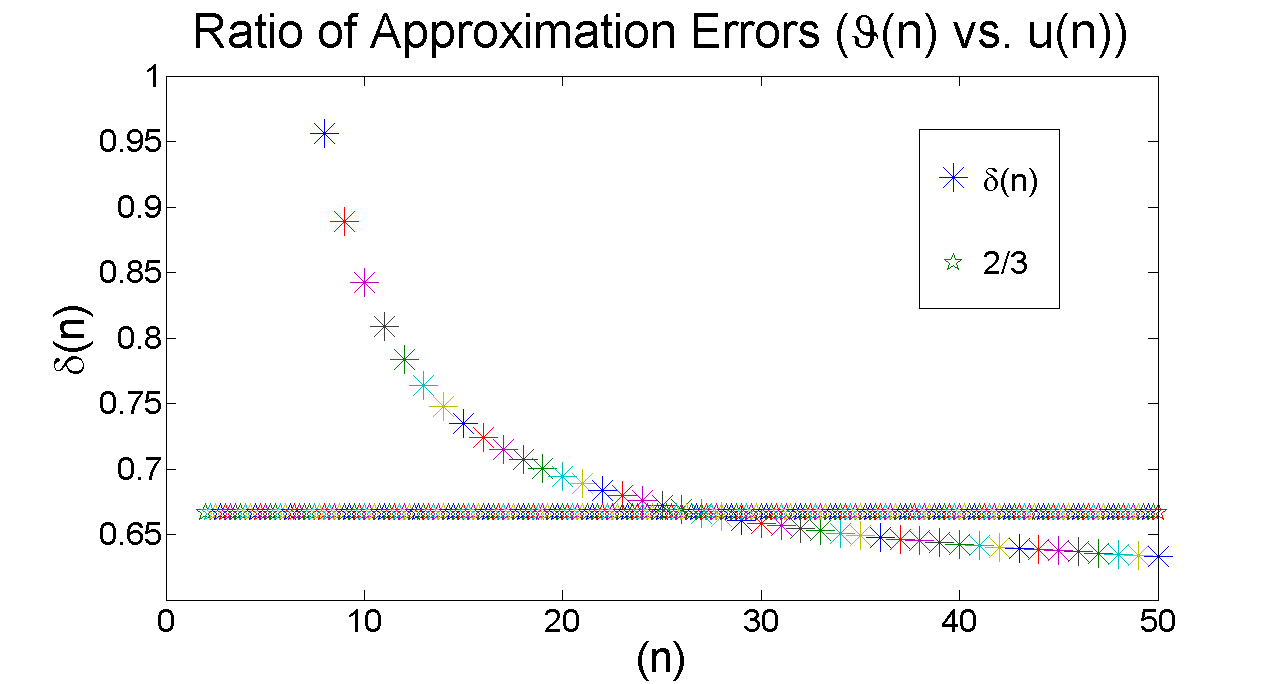} 
\caption{Ratio of the errors in approximation to $S_n$ associated with the estimates $\vartheta (n)$ and $u(n)$.  At $n=28$ the ratio $\delta(n)$ drops below $(2/3)$.}
\label{fig1}
\end{figure}

In this section we evaluate how the new estimate $(\ref{cat4})$ improves the approximation of the asymptotic limit $(\ref{cat3})$. Consider the ratio of the errors in approximation, \begin{equation}\label{rat}   \delta(n) = \frac{ S_n - \vartheta(n) }{ u(n) - S_n } \ .  \end{equation} The error in the approximation to $S_n$ by $\vartheta(n)$ is lower than that obtained from $u(n)$ by a factor of $1/\delta(n)$. Accordingly, values of $\delta(n)$ near zero imply $(\ref{cat4})$ is a significantly better estimate than $(\ref{cat3})$; note that $(\ref{main})$ implies $\delta(n) \in [0,1)$ for integers $n \geq 8$. In Fig.$\ref{fig1}$ the ratio $\delta(n)$ is plotted for integers $n \in [8,9,\ldots,50]$. At $n = 28$ the ratio $\delta(n)$ drops below $(2/3)$, which is plotted as a horizontal line.

To put our results in context, we compare $\delta(n)$ with a similar measure used in \cite{dut86}. The upper bound $(\ref{ub})$ was proven in \cite{dut86} to approximate $C_k$ at least 3 times as well as the previously established estimate $\upsilon(k) \doteq 4^k/((k+1) \sqrt{\pi k})$ \cite{com70}, $$ ( \nu(k) - C_k ) \leq \frac{1}{3} ( \upsilon(k) - C_k)  \ . $$   From Fig.$\ref{fig1}$ we find $\delta(n)$ drops below (2/3) for $n \geq 28$, thus our estimate $\vartheta_n$ improves the established estimate $u_n$ comparably to the improvement of $\nu_k$ over $\upsilon_k$ that was proven in \cite{dut86}.

\begin{figure}[htb] \center
\includegraphics[width=1 \linewidth]{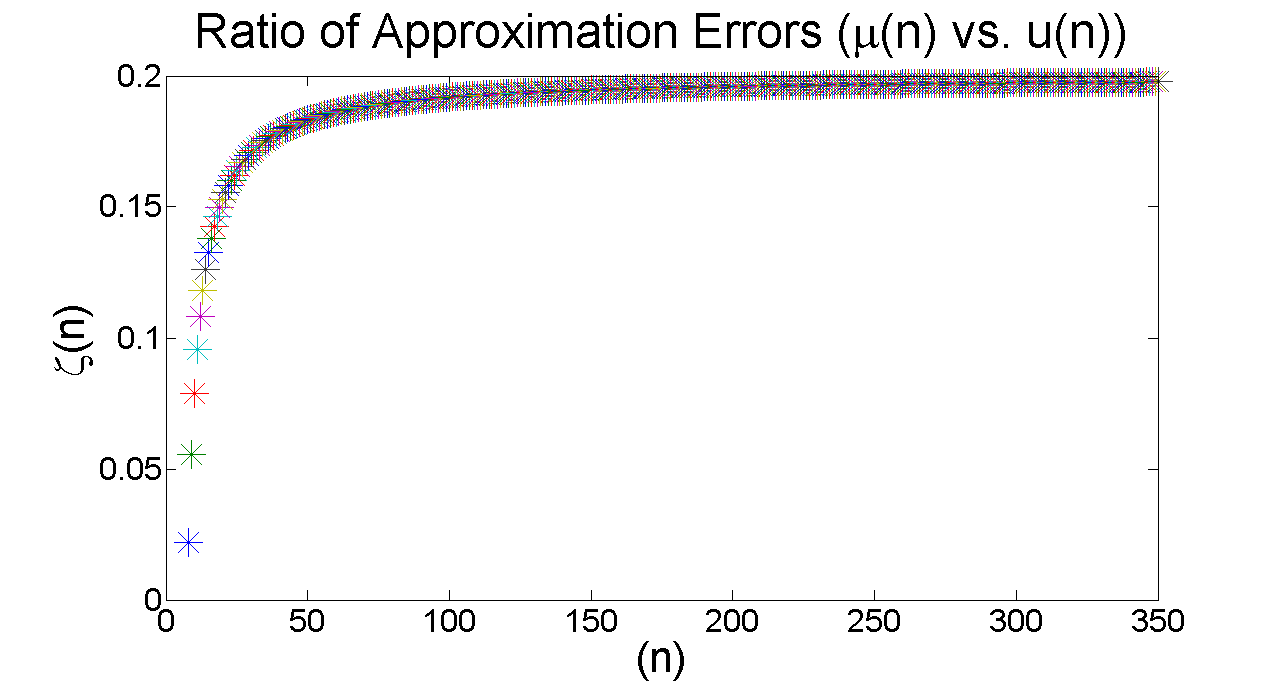} 
\caption{Ratio of the errors in approximation to $S_n$ associated with the estimates $\mu (n)$ and $u(n)$. The ratio $\zeta (n)$ approaches $(1/5)$ from below as $n \rightarrow \infty$.}
\label{fig2}
\end{figure}

In \cite{dut86} both a lower and upper bound on the $k^{th}$ Catalan number $C_k$ was established. In the numerical results presented in\cite{dut86} it was found that the average of the lower and upper bound significantly improved the approximation of $C_k$. This motivates us to consider taking the average of the lower bound $\vartheta(n)$ (cf.$(\ref{cat4})$) and upper bound $u(n)$ (cf.$(\ref{cat3})$) as an approximation to $S_n$. Define $\mu(n) \doteq  \frac{1}{2} \big( \vartheta(n) + u(n) \big)$. In Fig.$\ref{fig2}$ we plot the ratio of errors in approximation, $$  \zeta(n) = \frac{ \mu(n) - S_n }{ u(n) - S_n } \ . $$

From Fig.$\ref{fig2}$ we find that the ratio $\zeta (n)$ approaches $(1/5)$ from below as $n \rightarrow \infty$. In Fig.$\ref{fig1}$, the ratio $\delta(n)$ (cf.$(\ref{rat})$) approaches a value larger than $(3/5)$ from above as $n$ grows, thus we have $\zeta(n) < (1/3) \delta(n)$ and, consequently, the average $\frac{1}{2} \big( \vartheta(n) + u(n) \big)$ improves the estimate of $\vartheta(n)$ by a factor greater than $3$.

\begin{figure}[htb] \center
\includegraphics[width=1 \linewidth]{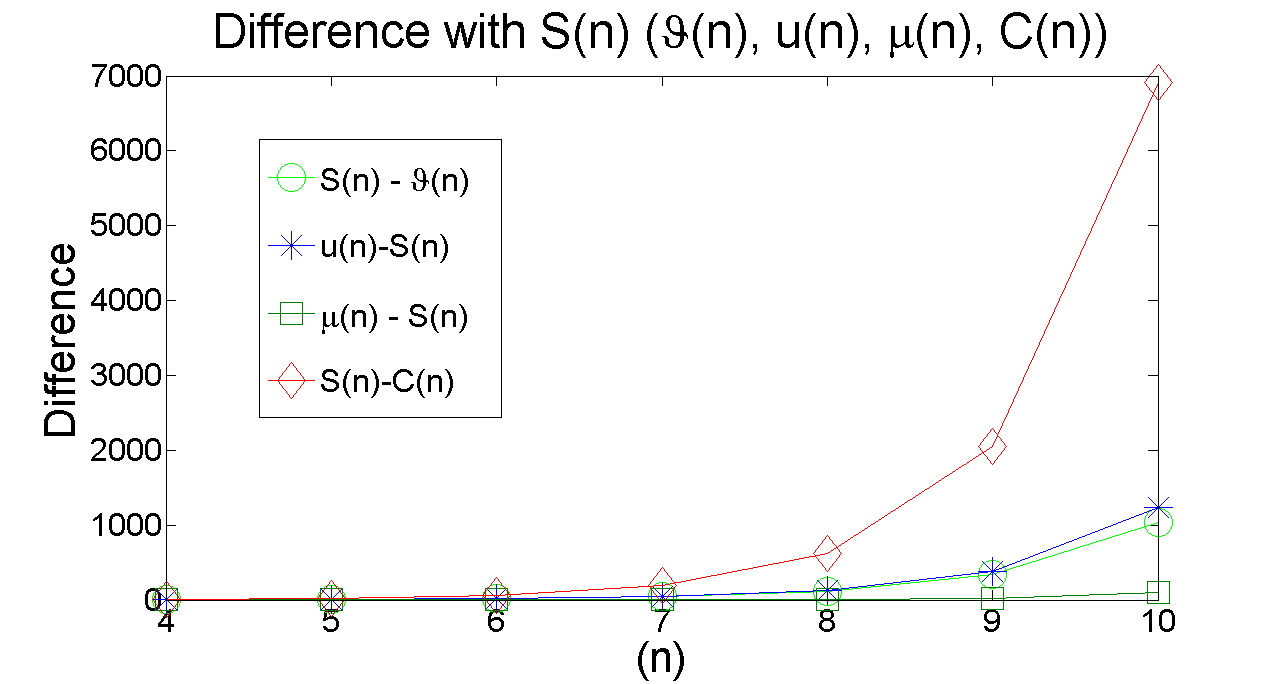} 
\caption{Difference between $S_n$ and $\{  \vartheta(n), u(n), \mu(n),C_n  \}$.}
\label{fig3}
\end{figure}

In Fig.$\ref{fig3}$, we plot the difference between $S_n$ and the estimates $\{ \vartheta(n), u(n), \mu(n), C_n  \}$. Clearly $C_n$ is significantly smaller than $S_n$, whereas both the upper bound $u(n)$ and lower bound $\vartheta(n)$ provide relatively similar approximations to $S_n$, albeit $\vartheta(n)$ remains the better estimate. The average $\mu(n)$ is empirically shown to be the best estimator of $S_n$ among the set $\{ \vartheta(n), u(n), \mu(n), C_n  \}$, as was suggested by Fig.$\ref{fig1},\ref{fig2}$. Note that $(\ref{main})$ implies $\mu(n)$ is also an upper bound to $S_n$.

% Empirically then we have $\mu(n)$ is at least $(9/2)$ times as good an upper bound to $S_n$ than the previously established estimate, which by $(\ref{main})$ is also an upper bound, $u(n)$.

\section{Future Work}\label{sec:fut}

It is possible to obtain an asymptotic approximation to $S_n$ that is arbitrarily tight by utilizing the following recurrence relation proposed in \cite{mather}, \begin{equation}\label{rec} (n+1) S_n + (1-5n)S_{n-1} = 2(1-2n)S_{n-2} \ .  \end{equation}  Specifically, by substituting the sum $S_n$ in $(\ref{rec})$ with a $p^{th}$ degree polynomial in $(1/n)$ and upper bound $(\ref{cat3})$ coefficient, \begin{equation}\label{t0} S_n \approx u(n) \Bigg( \sum_{r = 0}^p  \Big( \frac{c_r}{n} \Big)^r  \Bigg)  \ , \end{equation} or even more accurately by a $p^{th}$ degree polynomial in $(1/n)$ and lower bound $(\ref{cat4})$ coefficient, \begin{equation}\label{t1}  S_n \approx \vartheta(n) \Bigg( \sum_{r = 0}^p  \Big( \frac{m_r}{n} \Big)^r  \Bigg) \end{equation} we can iteratively solve for $\{c_r \ : \ r \in \{1,2,\ldots, p\} \}$ (resp. $m_r$) and obtain an asymptotic estimate for $S_n$ that becomes arbitrarily tight as $p$ approaches infinity. For $p =0$, Fig.$\ref{fig1}$ illustrates that $(\ref{t1})$ yeilds an estimate of $S_n$ that is at least $(3/2)$ times tighter than that of $(\ref{t0})$. For $p=1$ it can be shown that this ratio increases to $2$. Such approximations require an increasing number of computations and thus do not benefit from the relative simplicity of $(\ref{cat3})$ and $(\ref{cat4})$. An interesting project might consider how increasing the value of $p$ in $(\ref{t0})-(\ref{t1})$ will affect the error ratio between the two estimates, particularly in regard to the extra computational costs.

\section{Conclusion}\label{sec:con}

We have proven upper and lower bounds on the sum of Catalan numbers, $S_n$, where previously only an asymptotic limit had been proposed \cite{oeis}. The lower bound was proven to be a better approximation to $S_n$ than the upper bound, and empirical evidence shows this improvement is at least by a factor of $(3/2)$. The improvement of the lower bound over the previously established upper bound was shown to be comparable to the improvement in approximation of the $k^{th}$ Catalan number that was presented in \cite{dut86}. Motivated by the results presented in \cite{dut86}, the average of these bounds was considered.  The average proved to be an upper bound to $S_n$, and, empirically, it was found that the average provided a significant improvement on both the upper and lower approximations. Specifically, the average improved the lower bound by a factor greater than $3$, and improved the upper bound by a factor greater than $(9/2)$.

% \begin{equation}\label{com1} {2n  \choose n } = \frac{2 (2n - 1)}{n} {2n - 2 \choose n - 1}  \end{equation}  

\end{document}